\documentclass[10pt]{article}

 \usepackage{graphicx}

\usepackage[english,francais]{babel}

\usepackage{amssymb}
\usepackage{amsmath}

\newtheorem{theorem}{Theorem}[section]

\newtheorem{e-proposition}[theorem]{Proposition}

\newtheorem{e-definition}[theorem]{Definition\rm}

\def\og{\leavevmode\raise.3ex\hbox{$\scriptscriptstyle\langle\!\langle$~}}
\def\fg{\leavevmode\raise.3ex\hbox{~$\!\scriptscriptstyle\,\rangle\!\rangle$}}

\usepackage[latin1]{inputenc}

\graphicspath{
{./Figures/}
}

\newcommand{\stokeshat}{{\hat A}_{Stokes}}

\newcommand{\R}{{\bf R}}

\def \be{\begin{equation}}
\def \ee{\end{equation}}
\def \bea{\begin{eqnarray}}
\def \eea{\end{eqnarray}}
\def \bean{\begin{eqnarray*}}
\def \eean{\end{eqnarray*}}

\begin{document}

\bibliographystyle{alpha}

\title{New constructions of domain decomposition methods for systems of PDEs}

\author{   
V.~Dolean\footnote{Univ. d'\'Evry and CMAP, Ecole Polytechnique, 91128  
Palaiseau, France,  dolean@cmap.polytechnique.fr},~  
F.~Nataf\footnote{CMAP, CNRS UMR 7641, Ecole Polytechnique, 91128 Palaiseau Cedex, France, nataf@cmap.polytechnique.fr},~
G. Rapin\footnote{Math. Dep., NAM, University of Göttingen, D-37083, Germany,grapin@math.uni-goettingen.de }
}  

\maketitle

\begin{abstract}
We propose new domain decomposition methods for systems of partial differential equations in two and three dimensions. The algorithms are derived with the help of the Smith factorization of the operator. This could also be validated by numerical experiments. 
\end{abstract}

\section{Introduction}
Neumann-Neumann \cite{Bourgat:1989:VFA} or FETI type algorithms are very popular domain decompositions methods. They are currently used for very large scale computations, see for example \cite{DDO:2003:P15} and references therein. These methods are very well understood for symmetric definite positive scalar equations. For nonsymmetric problems and systems of equations many questions are still open. We propose in this note, a systematic construction of related algorithms for systems of partial differential equations (PDE). The approach is based on the Smith factorization of the system of PDEs. First, we explain the derivation of the domain decomposition method for the Stokes system. Then the application to the Oseen system \cite{Nataf:2005:NDS} and the Euler equations \cite{Dolean:2005:NDE} is briefly discussed. 
\section{Analysis of the Stokes system via Smith factorization}\label{sec:smith}
Let  $\nu>0$, we write the 2D Stokes equations as:
\be\label{eq:stokes}
A_{Stokes}\ \left( \begin{array}{c} u\\ v\\ p \end{array} \right)  = 
\left(\begin{array}{ccc} -\nu \Delta  & 0  &  - \partial_x \\
                         0 & -\nu \Delta  &  - \partial_y \\
                         \partial_x  & \partial_y & 0  
\end{array} \right)
\ \left( \begin{array}{c} u\\ v\\ p \end{array} \right)
= \ \left( \begin{array}{c} f_u\\ f_v\\ 0 \end{array} \right).
\ee
We first recall the definition of the Smith factorization of a matrix with polynomial entries and apply it to systems of PDEs: 
\begin{theorem}
Let $n$ be an integer and $A$ an invertible $n\times n$ matrix with polynomial entries with respect to the variable $\lambda$: $A=(a_{ij}(\lambda))_{1\le i,j\le n}$. \\
Then, there exist  matrices  $E$, $D$ and $F$ with polynomial entries satisfying the following properties:
det($E$)=det($F$)=1, $D$ is a diagonal matrix and $A=EDF$.
\end{theorem}
More details can be found in \cite{Wloka_Rowley_Lawruk:95}. We first take formally the Fourier transform of system (\ref{eq:stokes}) with respect to $y$ (dual variable is $k$). We keep the partial derivatives in $x$ since in the sequel we shall consider a model probem where the interface between the subdomain is orthogonal to the $x$ direction. We note
\begin{equation}\label{eq:stokeshat}
\stokeshat = 
\left(\begin{array}{ccc} -\nu (\partial_{xx}-k^2) & 0  &  - \partial_x \\
                         0 & -\nu (\partial_{xx}-k^2)  &  - ik  \\
                         \partial_x  & ik & 0  
\end{array} \right).
\end{equation}
We can perform the Smith factorization of  $\stokeshat$ by considering it as a matrix with polynomials in $\partial_x$ entries. We have  
\begin{equation}\label{eq:smithfactorization}
\stokeshat=EDF
\end{equation}
 where
$D_{11}=D_{22}=1 \text{ and }D_{33}= -\nu \hat{\Delta}^2$ with $\hat{\Delta} = \partial_{xx}-k^2$. One should note that a stream function formulation gives the same differential equation for the stream function. In the same way, the three-dimensional case can be characterized. In this case, the diagonal matrix $D_{3D}$ is a four by four matrix whose entries are:
${D_{3D}}_{,11}={D_{3D}}_{,22}=1$, ${D_{3D}}_{,33}= -\nu \hat{\Delta}$ and  ${D_{3D}}_{,44}= -\nu \hat{\Delta}^2$. This suggests that the derivation of a DDM for the bi-Laplacian is a key ingredient for a DDM for the Stokes system.

\section{A Domain Decomposition Method for the bi-Lplacian}
Let $\Omega$ be an open subset of $\R^2$ and $\Gamma=\bar\Omega\cap\ \{x=0\}$ be a symmetry axis of $\Omega$. For simplicity, in the note we assume homogeneous Dirichlet conditions on the boundary $\partial\Omega$. The domain $\Omega$ is decomposed into $\Omega_1=\Omega\cap  \{x<0\}$ and $\Omega_2=\Omega\cap  \{x>0\}$. We consider the following algorithm:\\
Starting with an initial guess satisfying $w_1^0=w_2^0$ and $\Delta(w_1^0)=\Delta(w_2^0)$ on $\Gamma$, the correction steps are expressed as follows for $i=1,2$:
\begin{eqnarray}
-\nu\Delta^2(\bar w_i^{k+1})&=&0 \text{ in } \Omega_i  \label{eq:corredp}\\
\partial_{n_i}\Delta \bar w_i^{k+1}&=&-(\partial_{n_1}\Delta w_1^k+\partial_{n_2}\Delta w_2^k)/2 \label{eq:corrbc}\\
\text{ and } \partial_{n_i} \bar w_i^{k+1}&=&-(\partial_{n_1} w_1^k+\partial_{n_2} w_2^k)/2 \text{ on }\Gamma \label{eq:corrbcbis}
\end{eqnarray}
followed by an update step:
\begin{eqnarray}
-\nu\Delta^2( w_i^{k+1})&=&g \text{ in } \Omega_i  \label{eq:updateedp}\\
\Delta w_i^{k+1}&=&\Delta w_i^{n}+(\Delta \bar w_1^{k+1}+\Delta \bar w_2^{k+1})/2 \label{eq:updatebc}\\
\text{ and }  w_i^{k+1}&=&w_i^k+(\bar w_1^{k+1}+\bar w_2^{k+1})/2  \text{ on }\Gamma \label{eq:updatebcbis}
\end{eqnarray}
By symmetry arguments, we have converge in two steps to the solution of $-\nu\Delta^2(w)=g$ in $\Omega$. 
\section{The algorithm for the Stokes system}\label{sec:newalgo}
Thanks to the Smith factorization (\ref{eq:smithfactorization}), it is possible to translate the above algorithm for the bi-Laplacian operator into an algorithm for the Stokes system. It suffices to replace equations (\ref{eq:corredp}),(\ref{eq:updateedp}) by  the Stokes equations and in the interface conditions (\ref{eq:corrbc}),(\ref{eq:corrbcbis}),(\ref{eq:updatebc}),(\ref{eq:updatebcbis}) $w$ by the last component of $F(u,v,p)^T$. In order to write the resulting algorithm in an intrinsic form, we introduce the stress $\vec\sigma(\vec u,p)$ on the interface for a velocity $\vec u=(u,v)$ and a pressure $p$. For any vector $\vec u$ its normal (resp. tangential) component on the interface is $\vec u_n$ (resp. $\vec u_\tau$). We denote $\vec\sigma_n$ and $\vec\sigma_\tau$ the normal and tangential parts of $\vec\sigma$, respectively. The new algorithm for the Stokes system for the same geometry as above reads:\\
Starting with an initial guess satisfying $\vec u_{1,\tau_1}^0=\vec u_{2,\tau_2}^0$ and $\vec\sigma_{1,n_1}^0=-\vec\sigma_{2,n_2}^0$ on $\Gamma$, the correction steps is expressed as follows for $i=1,2$:
\begin{eqnarray}
A_{Stokes}(\bar {\vec w}_i^{k+1},\bar q_i^{k+1})^T&=&0 \text{ in } \Omega_i  \label{eq:stokescorredp}\\
 \bar {\vec w}_{i,n_i}^{k+1}&=&(-1)^i(\vec u_{1,n_1}^k+ \vec u_{2,n_2}^k)/2 \\
\text{ and } 
\sigma_{\tau_i}(\vec {\bar w}_i^{k+1},\bar q_i^{k+1})&=&-(\sigma_{\tau_1}(\vec {\bar u}_1^k,\bar p_2^k)+\sigma_{\tau_2}(\vec {\bar u}_2^k,\bar p_2^k))/2  \text{ on }\Gamma \label{eq:stokescorrbc}
\end{eqnarray}
followed by an update step:
\begin{eqnarray}
A_{Stokes}(\vec u_i^{k+1},p_i^{k+1})^T&=&f \text{ in } \Omega_i  \label{eq:stokesupdateedp}\\
\vec u_{i,\tau_i}^{k+1}&=&\vec u_{i,\tau_i}^k+(\bar{\vec w}_{1,\tau_1}^{k+1}+\bar{\vec w}_{2,\tau_2}^{k+1})/2  \text{ on }\Gamma \\
\vec\sigma_{n_i}(\vec u_i^{k+1},p_i^{k+1})&=&\vec\sigma_{n_i}(\vec u_i^k,p_i^k)\nonumber\\
+ (-1)^{i-1}(\vec\sigma_{n_1}(\bar{\vec w}_1^{k+1},\bar q_1^{k+1})&-&\vec\sigma_{n_2}(\bar{\vec w}_2^{k+1},\bar q_2^{k+1}))/2  \text{ on }\Gamma \label{eq:stokesupdatebc}
\end{eqnarray}
The boundary conditions in the correction step involve the normal velocity and the tangential stress whereas in the update step they involve the tangential velocity and the normal stress. In 3D, the algorithm has the same definition. By construction, it converges in two steps. In the iterative version of the Neumann-Neumann algorithm for the Stokes system \cite{LeTallec:Patra:97}, \cite{Widlund:Pacvarino:02}, the boundary conditions of the correction step involve all the components of the stress whereas the update step involves all the components of the velocity. It can be shown that the convergence in two steps is then lost. More precisely, one obtains a convergence rate of $1/3$ in the case $\Omega=\R^2$, cf \cite{Nataf:2005:NDS}. 
\section{Algorithm for other systems of PDEs}
The derivation of the algorithm for the Stokes system is based on the use of the Smith factorization and on the existence of superconvergent algorithms for scalar PDEs. The same procedure can be performed for the Oseen equations \cite{Nataf:2005:NDS}. In this case, the diagonal form of the operator is the product of a convection-diffusion operator and of a Laplacian operator. Using (\cite{Achdou:1997:RPA})-\cite{ACHDOU:dd_p:00}) for the convection-diffusion, it is possible to derive an algorithm for the fourth order problem that converges in two steps. Translating this algorithm on the system, we obtain an algorithm converging in two steps for the Oseen system. The same work was done for the compressible Euler system \cite{Dolean:2005:NDE}. In this case, the diagonal form of the operator is a product of a convective operator with a Helmholtz convective wave equation. The Smith factorization has also been used to design PML for the time-dependent compressible Euler equations, \cite{Nataf:2005:NCP}. 
\section{Preliminary Numerical results for the Stokes system}
The domain $\Omega=(-A,B)\times (0,1)$ is decomposed into two subdomains $\Omega_1=(-A,0)\times (0,1)$ and  $\Omega_2=(0,B)\times (0,1)$. We compare the algorithm of \S~\ref{sec:newalgo} to the iterative version of the Neumann-Neumann algorithm. The stopping criteria is that the jumps of the normal derivative of the tangential component of the velocity are reduced by the factor $10^{-4}$. In table~\ref{tab:h} (left) $A=B=1$, we see that both algorithms are not sensitive with respect to the mesh size. Of course, due to the discrete approximation we cannot expect the optimal convergence in two steps. But we only need one more step to achieve the error bound. We have also varied the width of the subdomains, (middle table). As expected the convergence of the Neumann-Neumann method deteriorates. For large aspect ratios, the method diverges (-- in the table), since there exists an eigenvalue of the operator corresponding to the Richardson iteration with a modulus larger than $1$. But still in this case convergence can be enforced by its use as a preconditioner in Krylov method as it is usually the case. Our new algorithm seems to be  surprisingly robust with respect to the subdomain widths. For moderate variations we always need $3$ iterations steps. If we choose very thin subdomains, for instance $A=1$, $B=20$, the stopping criterion is achieved in only $7$ steps. In table~\ref{tab:h} (right), we have added a reaction term $c>0$ to the first two equations of the Stokes system. For instance $c$ may be the inverse of the time step in a time-dependent computation. We see that the new algorithm is fairly stable. 
\begin{table}
\begin{center}\hfill 
\begin{tabular}{|c||c|c|}
\hline
$h$ & new alg & N-N \\
\hline
0.02  & 3 & 10 \\
0.025 & 3 & 12 \\
0.05  & 3 & 11 \\
0.5   & 3 & 11 \\
0.1   & 3 & 11 \\
0.2   & 3 & 10 \\
\hline
\end{tabular}
\hfill
\begin{tabular}{|c||c|c|}
\hline
$B$ & new alg & N-N \\
\hline
1  & 3 & 11 \\
2 & 3 & 12 \\
3  & 3 & 11 \\
5   & 3 & 15 \\
10   & 3 & -- \\
20   & 7 & -- \\
\hline
\end{tabular}
\hfill
\begin{tabular}{|c||c|c|}
\hline
$c$ & new alg & N-N \\
\hline
0.001  & 3 & 11 \\
0.01   & 3 & 16 \\
0.1    & 3 & 19 \\
1      & 3 & 19 \\ 
10     & 3 & 16 \\
100    & 3 & 10\\
\hline
\end{tabular}
\hfill \\
\caption{Comparison between the new algorithm and the Neumann-Neumann algorithm (NN): Iteration counts for different mesh sizes (left), aspect ratio (middle) and different reaction terms (right)}  
\label{tab:h}
\end{center}
\end{table}

\bibliography{../ddm_stokes_3d/ddm_database,notesStokesO}

\end{document}